\documentclass[12pt, reqno]{amsart}
\usepackage{amsmath, amsthm, amscd, amsfonts, amssymb, graphicx, color}
\usepackage[bookmarksnumbered, colorlinks, plainpages]{hyperref}
\hypersetup{colorlinks=true,linkcolor=red, anchorcolor=green, citecolor=cyan, urlcolor=red, filecolor=magenta, pdftoolbar=true}

\newtheorem{theorem}{Theorem}[section]
\newtheorem{lemma}[theorem]{Lemma}

\newtheorem{corollary}[theorem]{Corollary}
\theoremstyle{definition}

\theoremstyle{remark}
\newtheorem{remark}[theorem]{Remark}
\numberwithin{equation}{section}

\newcommand{\Qcoh}{\hbox{\bf Qcoh}}
\newcommand{\QGr}{\hbox{\bf QGr}}

\begin{document}
\setcounter{page}{1}

\title[Rational elliptic curves]{Noncommutative geometry of  rational elliptic curves}

\author[I.Nikolaev]
{Igor ~V.~ Nikolaev$^1$}

\address{$^{1}$ 
Department of Mathematics and Computer Science, St.~John's University, 8000 Utopia Parkway,  
New York,  NY 11439, United States}

\email{\textcolor[rgb]{0.00,0.00,0.84}{igor.v.nikolaev@gmail.com}}



\date{Received:  February 12, 2017; Revised: yyyyyy; Accepted: zzzzzz.}

\begin{abstract}
We study an interplay between  operator algebras and 
geometry  of rational  elliptic curves.
Namely,  let $\mathcal{O}_B$ be  the Cuntz-Krieger algebra  given by  
square matrix  $B=(b-1, ~1, ~b-2, ~1)$,  where $b$ is an integer 
greater or equal to two.  It is proved,  that there exists a dense self-adjoint
sub-algebra of $\mathcal{O}_B$,  which is isomorphic  (modulo an ideal) to a 
twisted homogeneous   coordinate ring  of the rational elliptic curve 
$\mathcal{E}({\Bbb Q})=\{(x,y,z) \in {\Bbb P}^2({\Bbb C}) ~|~  y^2z=x(x-z)(x-{b-2\over b+2}z)\}$. 
\end{abstract}

\maketitle

\section{Introduction}
In 1950's due to the works of J.-P.~Serre and P.~Gabriel,  it was realized
that algebraic geometry can be recast in terms of noncommutative 
algebra;  we refer the reader to excellent survey by  [Stafford \& van~ den ~Bergh 2001]  \cite{StaVdb1}.
The following simple example  illustrates  the idea.     If $X$ is
a Hausdorff topological space and $C(X)$ is the $C^*$-algebra of continuous complex valued functions on $X$, 
then by the Gelfand Theorem the topology of $X$ is determined   by the commutative algebra $C(X)$.
This fact   can be written as   $K_0^{top}(X)\cong K_0^{alg}(C(X))$,
where $K_0^{top}$ and $K_0^{alg}$ are the topological and algebraic $K_0$-groups, respectively [Blackadar 1986]  \cite{B}.
 Consider now   the algebra $C(X)\otimes M_2({\Bbb C})$ consisting of 
  the two-by-two matrices  with entries in $C(X)$.
  Since the algebraic K-theory is stable  under the tensor products, one gets an isomorphism    
$K_0^{alg}(C(X))\cong K_0^{alg}(C(X)\otimes M_2({\Bbb C}))$
[Blackadar 1986]  \cite[\S 5]{B}. 
In other words,  topology of the space $X$ is determined by 
the tensor product $C(X)\otimes M_2({\Bbb C})$ which is 
 no longer a commutative algebra.
In the context of algebraic geometry,  one replaces  the space $X$
by a  projective variety $V$,  the algebra $C(X)$ by the coordinate ring of $V$,
the tensor product $C(X)\otimes M_2({\Bbb C})$ by a twisted coordinate ring of $V$  and the group  $K^{top}(X)$
by a category of the quasi-coherent sheaves on $V$  [Stafford \& van~ den ~Bergh 2001]  \cite[p.173]{StaVdb1}.  
Below we  give a brief review  of this construction when $V$ is an elliptic curve;  we refer the reader to 
[Sklyanin 1982]  \cite{Skl1},   [Smith \& Stafford 1992]  \cite[pp.265-268 ]{SmiSta1} and 
[Stafford \& van~ den ~Bergh 2001]  \cite[p.197]{StaVdb1} for a detailed account.

Let $k$ be a field of $char~(k)\ne 2$. The  {\it Sklyanin algebra} 
$\mathfrak{S}_{\alpha,\beta,\gamma}(k)$
is a free $k$-algebra  on  four generators  $x_i$ and 
six  quadratic relations:
\begin{equation}\label{eq1}
\left\{
\begin{array}{ccc}
x_1x_2-x_2x_1 &=& \alpha(x_3x_4+x_4x_3),\\
x_1x_2+x_2x_1 &=& x_3x_4-x_4x_3,\\
x_1x_3-x_3x_1 &=& \beta(x_4x_2+x_2x_4),\\
x_1x_3+x_3x_1 &=& x_4x_2-x_2x_4,\\
x_1x_4-x_4x_1 &=& \gamma(x_2x_3+x_3x_2),\\ 
x_1x_4+x_4x_1 &=& x_2x_3-x_3x_2,
\end{array}
\right.
\end{equation}
where $\alpha,\beta,\gamma\in k$ and 
$\alpha+\beta+\gamma+\alpha\beta\gamma=0$.
If $\alpha\not\in \{0;\pm 1\}$  then  algebra $\mathfrak{S}_{\alpha,\beta,\gamma}(k)$  defines a
non-singular  elliptic  curve $\mathcal{E}\subset {\Bbb P}^3(k)$ given by an intersection of the quadrics
$u^2+v^2+w^2+z^2={1-\alpha\over 1+\beta}v^2+{1+\alpha\over 1-\gamma}w^2+z^2=0$
together with an automorphism $\sigma: \mathcal{E}\to \mathcal{E}$.
We shall use the following isomorphism (see [Sklyanin 1982]  \cite{Skl1} and [Smith \& Stafford 1992]  \cite{SmiSta1}):
\begin{equation}\label{eq2}
\QGr~ (\mathfrak{S}_{\alpha,\beta,\gamma}(k)~/~\Omega)\cong  \Qcoh~(\mathcal{E}),
\end{equation}
where $\QGr$ is a category of the quotient graded modules over 
the algebra $\mathfrak{S}_{\alpha,\beta,\gamma}(k)$ modulo torsion, $\Qcoh$ a category 
of the quasi-coherent
sheaves on $\mathcal{E}$ and $\Omega\subset \mathfrak{S}_{\alpha,\beta,\gamma}(k)$ 
a two-sided ideal generated by the  central  elements
$\Omega_1 = -x_1^2+x_2^2+x_3^2+x_4^2$ and 
$\Omega_2 = x_2^2+{1+\beta\over 1-\gamma}x_3^2+{1-\beta\over 1+\alpha}x_4^2$
[Smith \& Stafford 1992]  \cite[p.276]{SmiSta1}.  The quotient  of the Sklyanin algebra by the ideal $\Omega$  
is called a {\it twisted   homogeneous coordinate ring} of the elliptic curve $\mathcal{E}$.

Let $A$ be a two-by-two matrix with non-negative
integer entries $a_{ij}$, such that every row and every column of $A$
is non-zero. 
The two-dimensional {\it Cuntz-Krieger algebra} $\mathcal{O}_A$
is  a $C^*$-algebra of bounded linear operators on a Hilbert space
$\mathcal{H}$  generated by the partial isometries $s_1$ and $s_2$,
and  relations:
\begin{equation}\label{eq3}
\left\{
\begin{array}{ccc}
s_1^*s_1 &=& a_{11}s_1s_1^*+a_{12}s_2s_2^*,\\
s_2^*s_2 &=& a_{21}s_1s_1^*+a_{22}s_2s_2^*,\\
Id &=& s_1s_1^*+s_2s_2^*,
\end{array}
\right.
\end{equation}
where $Id$ is the identity operator on $\mathcal{H}$.  
Occasionally,   the algebra $\mathcal{O}_A$ will be written as 
$\mathcal{O}_{a_{11}, a_{12}, a_{21}, a_{22}}$. 
If one defines $x_1=s_1, x_2=s_1^*,  x_3=s_2$ and  $x_4=s_2^*$, then it is 
easy to see, that   $\mathcal{O}_A$ contains a dense sub-algebra $\mathcal{O}_A^0$,
which is a free ${\Bbb C}$-algebra on four generators $x_i$ and three quadratic relations:
\begin{equation}\label{eq4}
\left\{
\begin{array}{ccc}
x_2x_1 &=& a_{11}x_1x_2+a_{12}x_3x_4,\\
x_4x_3 &=& a_{21}x_1x_2+a_{22}x_3x_4,\\
1 &=& x_1x_2+x_3x_4,
\end{array}
\right.
\end{equation}
and an involution acting by  the formula:
\begin{equation}\label{eq5}
x_1^*=x_2,  \qquad  x_3^*=x_4. 
\end{equation}
Notice,  that equations (\ref{eq4}) are invariant  under   this   involution.

It  is known,  that 
ideal (\ref{eq1}) is stable under  involution (\ref{eq5}), if
and only if,   $\bar\alpha=\alpha, \beta=1$ and $\gamma=-1$ (lemma \ref{lm1});
the involution turns the  Sklyanin algebra   $\mathfrak{S}_{\alpha, 1, -1}({\Bbb C})$ 
in a $\ast$-algebra (i.e a self-adjoint algebra).  Denote by $\mathcal{I}_0$ a  
non-homogeneous two-sided ideal
of  $\mathfrak{S}_{\alpha, 1, -1}({\Bbb C})$
generated by relation $x_1x_2+x_3x_4=1$.  
Let $\mathcal{J}_0$ be a two-sided  ideal of $\mathcal{O}_A^0$  
generated by four relations $x_4x_2-x_1x_3=x_3x_1+x_2x_4=x_4x_1-x_2x_3=
x_3x_2+x_1x_4=0$. 
 The following theorem  and corollary  describe a family of
 the Cuntz-Krieger algebras  which are  twisted homogeneous coordinate rings
of   the rational  elliptic curves. 
\begin{theorem}\label{thm1}
For every integer  $b\ge 2$,  there exists a $\ast$-isomorphism:
\begin{equation}\label{eq6}
\mathfrak{S}_{{b-2\over b+2}, ~1, ~-1}({\Bbb C})~/  ~\mathcal{I}_0 ~\cong
~\mathcal{O}^0_B~/ ~\mathcal{J}_0, \qquad \hbox{\it where} 
~B=\left(\small
\begin{matrix}
b-1 & 1\cr b-2 & 1
\end{matrix}
\right).
\end{equation}
\end{theorem}
\begin{corollary}\label{cor1}
For every integer  $b\ge 2$,   there exists 
a dense self-adjoint   sub-algebra of the Cuntz-Krieger algebra
$\mathcal{O}_B$ isomorphic  modulo the  ideal $\mathcal{I}_0$ to  
the  twisted homogeneous  coordinate ring of the rational elliptic curve
\begin{equation}  
\mathcal{E}_b({\Bbb Q})=\left\{(x,y,z) \in {\Bbb P}^2({\Bbb C}) ~|~  y^2z=x(x-z)(x-{b-2\over b+2}z)\right\}.
\end{equation}
\end{corollary}
\begin{remark}
There exists a canonical isomorphism 
\begin{equation}  
\mathcal{O}_B\otimes\mathcal{K}\cong {\Bbb A}_B\rtimes_{\sigma}{\Bbb Z},
\end{equation}
where ${\Bbb A}_B$ is an {\it AF-algebra} with the incidence matrix $B$ 
introduced by [Effros and Shen 1980]  \cite{EfSh1},  $\sigma$ is the shift
automorphism of  ${\Bbb A}_B$ and  $\mathcal{K}$ is the $C^*$-algebra of compact operators
 [Blackadar 1986]  \cite[Exercise 10.11.9]{B}.  
Thus  the algebra ${\Bbb A}_B$ is an analog of  the coordinate ring of the  $\mathcal{E}_b({\Bbb Q})$. 
This  observation can be used to  calculate  traces 
of the Frobenius endomorphisms  in terms of the 
algebra ${\Bbb A}_B$  \cite{Nik1}.  
\end{remark}
Our note is organized as follows. Theorem \ref{thm1} is proved in Section 2.
The proof of corollary \ref{cor1} can be found in Section 3.  
All preliminary facts  have been introduced in Section 1;
we refer the reader to  [Cuntz \& Krieger 1980]   \cite{CuKr1} and 
[Stafford \& van~ den ~Bergh 2001]  \cite{StaVdb1} for the details.

\section{Proof of theorem \ref{thm1}}
We shall split the proof in a series of  lemmas.
\begin{lemma}\label{lm1}
The  ideal of free algebra ${\Bbb C}\langle x_1,x_2,x_3,x_4\rangle$ 
generated by   equations (\ref{eq1})  is  stable  under   involution (\ref{eq5}), 
if and only if,  $\bar\alpha=\alpha, \beta=1$ and $\gamma=-1$. 
\end{lemma}
\begin{proof}
(i) Let us consider the first two equations (\ref{eq1});
this pair is invariant under involution (\ref{eq5}).  Indeed, by the rules of composition
for an involution
\begin{equation}\label{eq7}
\left\{
\begin{array}{ccc}
(x_1x_2)^* &= x_2^*x_1^* &= x_1x_2,\\
(x_2x_1)^* &= x_1^*x_2^* &= x_2x_1,\\
(x_3x_4)^* &= x_4^*x_3^* &= x_3x_4,\\
(x_4x_3)^* &= x_3^*x_4^* &= x_4x_3.
\end{array}
\right.
\end{equation}
Since $\alpha^*=\bar\alpha=\alpha$,  the first two equation (\ref{eq1}) 
remain invariant under involution (\ref{eq5}).

\smallskip
(ii) Let us consider the middle pair of equations (\ref{eq1});
 by the rules of composition for an involution
\begin{equation}\label{eq8}
\left\{
\begin{array}{ccc}
(x_1x_3)^* &= x_3^*x_1^* &= x_4x_2,\\
(x_3x_1)^* &= x_1^*x_3^* &= x_2x_4,\\
(x_2x_4)^* &= x_4^*x_2^* &= x_3x_1,\\
(x_4x_2)^* &= x_2^*x_4^* &= x_1x_3.
\end{array}
\right.
\end{equation}
One can apply the involution to the first equation 
$x_1x_3-x_3x_1=\beta(x_4x_2+x_2x_4)$; then one gets
$x_4x_2-x_2x_4=\bar\beta(x_1x_3+x_3x_1)$.  But the second
equation says that $x_1x_3+x_3x_1=x_4x_2-x_2x_4$;  the last  
two equations are compatible if and only if $\bar\beta=1$. 
Thus, $\beta=1$. 

The second equation in involution writes as $x_4x_2+x_2x_4=x_1x_3-x_3x_1$;
the last equation  coincides with  the first equation for $\beta=1$.   

Therefore,  $\beta=1$ is necessary and sufficient for invariance of the middle pair of
equations (\ref{eq1}) with respect to  involution (\ref{eq5}).

\smallskip
(iii) Let us consider the last pair of equations (\ref{eq1});
 by the rules of composition for an involution
\begin{equation}\label{eq9}
\left\{
\begin{array}{ccc}
(x_1x_4)^* &= x_4^*x_1^* &= x_3x_2,\\
(x_4x_1)^* &= x_1^*x_4^* &= x_2x_3,\\
(x_2x_3)^* &= x_3^*x_2^* &= x_4x_1,\\
(x_3x_2)^* &= x_2^*x_3^* &= x_1x_4.
\end{array}
\right.
\end{equation}
One can apply the involution to the first equation 
$x_1x_4-x_4x_1 = \gamma(x_2x_3+x_3x_2)$; then one gets
$x_3x_2-x_2x_3=\bar\gamma(x_4x_1+x_1x_4)$.  But the second
equation says that $x_1x_4+x_4x_1 = x_2x_3-x_3x_2$;  the last  
two equations are compatible if and only if $\bar\gamma=-1$. 
Thus, $\gamma=-1$. 

The second equation in involution writes as $x_3x_2+x_2x_3=x_4x_1-x_1x_4$;
the last equation  coincides with  the first equation for $\gamma=-1$.   

Therefore,  $\gamma=-1$ is necessary and sufficient for invariance of the last pair of
equations (\ref{eq1}) with respect to  involution (\ref{eq5}).

\smallskip
(iv) It remains to verify that condition $\alpha+\beta+\gamma+\alpha\beta\gamma=0$
is satisfied by $\beta=1$ and $\gamma=-1$ for any $\alpha\in k$.
Lemma \ref{lm1} follows.

\end{proof}

\begin{lemma}\label{lm2}
Whenever $\alpha\ne 1$ there exists an invertible linear 
transformation with  rational coefficients which brings  
the system of equations
\begin{equation}\label{eq10}
\left\{
\begin{array}{ccc}
x_1x_2-x_2x_1 &=& \alpha(x_3x_4+x_4x_3),\\
x_1x_2+x_2x_1 &=& x_3x_4-x_4x_3
\end{array}
\right.
\end{equation}
to the form
\begin{equation}\label{eq11}
\left\{
\begin{array}{ccc}
x_2x_1 &=& (b-1)x_1x_2+x_3x_4,\\
x_4x_3 &=& (b-2)x_1x_2+ x_3x_4, 
\end{array}
\right.
\end{equation}
where $\alpha={b-2\over b+2}$. 
\end{lemma}
\begin{proof}
(i) Let us isolate $x_2x_1$ and $x_4x_3$ in equations (\ref{eq10});
for that,  we shall write (\ref{eq10}) in the form
\begin{equation}\label{eq12}
\left\{
\begin{array}{ccc}
x_2x_1+\alpha x_4x_3 &=& x_1x_2-\alpha x_3x_4,\\
x_2x_1+x_4x_3 &=& -x_1x_2+x_3x_4.
\end{array}
\right.
\end{equation}
Consider (\ref{eq12}) as  a linear system of equations relatively
$x_2x_1$ and $x_4x_3$;   since $\alpha\ne 1$,  it  has a unique
solution
\begin{equation}\label{eq13}
\left\{
\begin{array}{ccccc}
x_2x_1 &=& {1\over 1-\alpha} \left|

\begin{matrix}
x_1x_2-\alpha x_3x_4 & \alpha\cr
                                                                  -x_1x_2+x_3x_4 &  1
\end{matrix}                                                                  
                                                                  \right|
&=&                                                                  
{1+\alpha\over 1-\alpha}x_1x_2-{2\alpha\over 1-\alpha}x_3x_4,\\
&&&&\\
x_4x_3 &=& {1\over 1-\alpha} \left|
\begin{matrix}
1 & x_1x_2-\alpha x_3x_4 & \cr
                                                                   1 & -x_1x_2+x_3x_4
\end{matrix}                                                                   
                                                                   \right|
&=&
{-2\over 1-\alpha} x_1x_2+{1+\alpha\over 1-\alpha}x_3x_4.
\end{array}
\right.
\end{equation}

\medskip
(ii)  Let us substitute $\alpha={b-2\over b+2}$ in (\ref{eq13});  then one arrives 
at the following system of equations  given in the matrix form
\begin{equation}\label{eq14}
\left(
\begin{matrix}
x_2x_1\cr x_4x_3
\end{matrix}
\right)= 
\left(
\begin{matrix}
{b\over 2} & 1-{b\over 2}\cr -1-{b\over 2} & {b\over 2}
\end{matrix}
\right)
\left(
\begin{matrix}
x_1x_2\cr x_3x_4
\end{matrix}
\right). 
\end{equation}
It is verified directly,  that 
\begin{equation}\label{eq15}
\left(
\begin{matrix}
{1\over 2} & -{1\over 2}\cr 1 & 0
\end{matrix}
\right) 
\left(
\begin{matrix}
{b\over 2} & 1-{b\over 2}\cr -1-{b\over 2} & {b\over 2}
\end{matrix}
\right)
\left(
\begin{matrix}
0 & 1\cr -2 & 1
\end{matrix}
\right)=
\left(
\begin{matrix}
b-1 & 1\cr b-2 & 1
\end{matrix}
\right).
\end{equation}
In other words,     matrices  (\ref{eq11}) and  (\ref{eq14}) are  similar 
 in  the  matrix group $GL_2({\Bbb Q})$.   Lemma \ref{lm2} is proved. 
\end{proof}

\begin{lemma}\label{lm3}
If $b\ge 2$ is an integer, then there exists a $\ast$-isomorphism
\begin{equation}\label{eq16}
\mathfrak{S}_{{b-2\over b+2}, ~1, ~-1}({\Bbb C})~/  ~\mathcal{I}_0 ~\cong
~\mathcal{O}^0_{b-1,  ~1, ~b-2,   ~1}~/ ~\mathcal{J}_0;
\end{equation}
the isomorphism is given by identification of generators $x_i$ of
the respective algebras. 
\end{lemma}
{\it Proof.}
Since $b$ is  integer number,  one gets $\alpha={b-2\over b+2}$
is a rational number.  In particular, $\alpha$ is real,  i.e. $\bar\alpha=\alpha$;  thus, by lemma \ref{lm1},
algebra $\mathfrak{S}_{{b-2\over b+2}, ~1, ~-1}({\Bbb C})$ is a self-adjoint
Sklyanin algebra.

Recall,  that  ideal $\mathcal{I}_0$ is generated by relation 
\begin{equation}\label{eq17}
x_1x_2+x_3x_4=1,
\end{equation}
while ideal $\mathcal{J}_0$ is generated by the system of relations
\begin{equation}\label{eq18}
\left\{
\begin{array}{ccc}
x_1x_3 &=& x_4x_2,\\
x_3x_1 &=& -x_2x_4,\\
x_1x_4 &=& -x_3x_2,\\
x_4x_1 &=& x_2x_3.
\end{array}
\right.
\end{equation}
Notice,  that ideals $\mathcal{I}_0$ and $\mathcal{J}_0$ are stable
under   involution (\ref{eq5}).

By lemma \ref{lm2},  the first pair of equations in the system (\ref{eq1}) 
with $\alpha={b-2\over b+2}$ coincides with the first pair of equations in the 
system (\ref{eq4}) with $a_{11}=b-1, a_{12}=1, a_{21}=b-2$ and $a_{22}=1$. 
Thus, if  one complements system (\ref{eq1}) with equation (\ref{eq17}) and     
system (\ref{eq4}) with the system of equations (\ref{eq18}),  then one obtains 
the required $\ast$-isomorphism (\ref{eq16}).  
 Lemma \ref{lm3} is proved.
$\square$

\bigskip
Theorem \ref{thm1}  follows from lemma \ref{lm3}.

\bigskip
\begin{remark}\label{rmk1}
Ideals $\mathcal{I}_0$  and $\mathcal{J}_0$ do not depend 
on  ``modulus''   $b$ of the Sklyanin algebra 
$\mathfrak{S}_{{b-2\over b+2}, ~1, ~-1}({\Bbb C})$;  therefore,
algebra  $\mathcal{O}^0_{b-1,  ~1, ~b-2,   ~1}$  can be viewed
as a twisted  homogeneous coordinate ring of the elliptic curve $\mathcal{E}\subset {\Bbb P}^3({\Bbb C})$.
\end{remark}

\section{ Proof of corollary \ref{cor1}}
We shall split the proof in a series of lemmas, starting with the following 
elementary
\begin{lemma}\label{lm5}
If $\alpha$ is a real number different from $0$ and $1$, then the algebra 
$\mathfrak{S}_{\alpha, ~1, ~-1}({\Bbb C})~/~\Omega_0$
is the coordinate ring of a non-singular elliptic curve
$\mathcal{E}({\Bbb C})=\{(x,y,z) \in {\Bbb P}^2({\Bbb C}) ~|~  y^2z=x(x-z)(x-\alpha z)\}$. 
 \end{lemma}
\begin{proof}
Recall, that the Sklyanin algebra $\mathfrak{S}_{\alpha, ~1, ~-1}({\Bbb C})$ defines
an elliptic curve $\mathcal{E}\subset {\Bbb P}^3({\Bbb C})$ given by the intersection
of two quadrics ([Smith \& Stafford 1992]  \cite[p.267]{SmiSta1}):
\begin{equation}\label{eq19}
\left\{
\begin{array}{ccc}
(1-\alpha)v^2+(1+\alpha)w^2+2z^2 &=& 0,\\
u^2+v^2+w^2+z^2&=& 0. 
\end{array}
\right.
\end{equation}
We shall pass in (\ref{eq19}) from variables $(u,v,w,z)$ to
the new variables $(X,Y,Z,T)$ given by the formulas
\begin{equation}\label{eq20}
\left\{
\begin{array}{ccc}
u^2 &=& T^2,\\
v^2 &=& {1\over 2}Y^2-{1\over 2}Z^2-T^2,\\
w^2 &=& X^2+{1\over 2}Y^2-{1\over 2}Z^2-T^2,\\
z^2 &=& Z^2.
\end{array}
\right.
\end{equation}
Then equations (\ref{eq19}) take the form  
\begin{equation}\label{eq21}
\left\{
\begin{array}{ccc}
\alpha X^2+Z^2-T^2 &=& 0,\\
X^2+Y^2-T^2 &=& 0. 
\end{array}
\right.
\end{equation}
Let us consider another (polynomial) transformation
$(x,y)\mapsto (X,Y,Z,T)$ given by the formulas
\begin{equation}\label{eq22}
\left\{
\begin{array}{ccc}
X &=& -2y,\\
Y &=& x^2-1+\alpha,\\
Z &=& x^2+2(1-\alpha) x+1-\alpha,\\
T &=& x^2+2x+1-\alpha.
\end{array}
\right.
\end{equation}
Then both of  the equations (\ref{eq21}) give us the equation
$y^2=x(x+1)(x+1-\alpha)$,  which after a shift $x'=x+1$ 
takes the canonical form
\begin{equation}\label{eq23}
y^2=x(x-1)(x-\alpha).
\end{equation}
Using projective transformation $x={x'\over z'}$ $y={y'\over z'}$
in (\ref{eq23}), one gets the homogeneous equation of elliptic 
curve $\mathcal{E}$:
\begin{equation}\label{eq23bis}
y^2z=x(x-z)(x-\alpha z).
\end{equation}

Lemma \ref{lm5} follows.
\end{proof}

\begin{lemma}\label{lm6}
If  $b\ge 2$ is an integer,  then  there exists 
a dense self-adjoint  sub-algebra of the Cuntz-Krieger algebra
$\mathcal{O}_{b-1,  ~1, ~b-2,   ~1}$,   which is  related  (modulo  ideal $\mathcal{I}_0$) to  
a twisted homogeneous   coordinate ring of the rational elliptic curve  
$\mathcal{E}({\Bbb Q})=\{(x,y,z) \in {\Bbb P}^2({\Bbb C}) ~|~  y^2z=x(x-z)(x-{b-2\over b+2}z)\}$;
the curve is non-singular unless $b=2$. 
 \end{lemma}
\begin{proof}
If one  assumes  $\alpha={b-2\over b+2}$ in lemma \ref{lm5},
then: 
\begin{equation}\label{eq25}
\mathfrak{S}_{{b-2\over b+2}, ~1, ~-1}({\Bbb C})~/  ~\mathcal{I}_0 ~\cong
\mathcal{O}^0_{b-1,  ~1, ~b-2,   ~1}~/ ~\mathcal{J}_0.
\end{equation}
The RHS of  (\ref{eq25}) is a  sub-algebra
of the Cuntz-Krieger algebra  $\mathcal{O}_{b-1,  ~1, ~b-2,   ~1}$;
such an algebra is self-adjoint,  since the ideal  $\mathcal{J}_0$ is 
 invariant of the involution (\ref{eq5}).

The RHS of  (\ref{eq25}) is a dense  sub-algebra
of the Cuntz-Krieger algebra  $\mathcal{O}_{b-1,  ~1, ~b-2,   ~1}$, 
since $\mathcal{O}^0_{b-1,  ~1, ~b-2,   ~1}~/ ~\mathcal{J}_0$
 is dense in   $\mathcal{O}_{b-1,  ~1, ~b-2,   ~1}$.

On the other hand,  if $b\ne 2$,  the algebra $\mathcal{O}^0_{b-1,  ~1, ~b-2,   ~1}~/ ~\mathcal{J}_0$
is related  to the factor (by the ideal $\mathcal{I}_0$)  of the  coordinate ring   $\mathfrak{S}_{{b-2\over b+2}, ~1, ~-1}({\Bbb C})~/  ~\Omega$
of the non-singular curve $\mathcal{E}({\Bbb Q})=\{(x,y) \in {\Bbb P}^2({\Bbb C}) ~|~  y^2z=x(x-z)(x-{b-2\over b+2}z)\}$.
It is easy to see,  that the curve $\mathcal{E}({\Bbb Q})$  is singular if and only if $b=2$.  
Lemma \ref{lm6} is proved. 
\end{proof}

\bigskip
Corollary \ref{cor1}  follows from lemma \ref{lm6}.

\bibliographystyle{amsplain}


\end{document}